\documentclass[11pt,a4paper]{amsart}
\usepackage[utf8]{inputenc}
\usepackage{amsmath}
\usepackage{amsfonts}
\usepackage{amssymb}
\usepackage{graphicx}
\usepackage{booktabs}
\usepackage{tikz}
\usepackage[ruled,vlined]{algorithm2e}

\usepackage[left=2.2cm,right=2.2cm,top=2cm,bottom=2cm]{geometry}

\newcommand{\D}{{\mathop{}\!\mathrm{d}}} 
\newcommand{\R}{\mathbb{R}}
\newcommand{\Q}{\mathbb{Q}}

\newcommand{\N}{\mathbb{N}}

\newcommand{\E}{\mathbb{E}}

\newcommand{\cdotb}{\boldsymbol{\cdot}}
\numberwithin{equation}{section}  
\newtheorem{defn}{Definition}[section]

\newtheorem{rem}[defn]{Remark}

\newtheorem{prop}[defn]{Proposition}

\usepackage{hyperref}
\usepackage{xcolor}
\hypersetup{
    colorlinks,
    linkcolor={red!50!black},
    citecolor={blue!50!black},
    urlcolor={blue!80!black}
}

\title[A Multi-Marginal Convex Duality Theorem for Martingale Optimal Transport]{A Multi-Marginal C-Convex Duality Theorem for Martingale Optimal Transport}

\author[J. Sester]{Julian Sester}

\begin{document}

\maketitle

\begin{center}
\normalsize{\today} \\ \vspace{0.5cm}
\small\textit{
National University of Singapore, Department of Mathematics,\\ 21 Lower Kent Ridge Road, 119077.}                                                                                                                              
\end{center}

\begin{abstract} 
A convex duality result for martingale optimal transport problems with two marginals was established in \cite{beiglbock2013model}. In this paper we provide a generalization of this result to the multi-period setting. 
\end{abstract}

%
\maketitle

\section{Introduction}
Given some measurable function $c: \R^n \rightarrow \R$ and fixed marginal distributions\footnote{Here $\mathcal{M}_1(\R)$ denotes the set of probability measures on the reals equipped with the Borel sigma-algebra. Moreover, we implicitly assume that the marginals increase in convex order to ensure the non-emptiness of $\mathcal{M}(\mu_1,\dots,\mu_n)$, compare also \cite{strassen1965existence}.} $\mu_1, \dots, \mu_n \in \mathcal{M}_1(\R)$ the \emph{martingale optimal transport (MOT)} problem  consists in solving
\begin{equation}\label{eq_min_MOT}
\inf_{\Q\in \mathcal{M}(\mu_1,\dots,\mu_n)} \int_{\R^n} c(x_1,\dots,x_n) ~\D \Q(x_1,\dots,x_n),
\end{equation}
where { \begin{align*}
\mathcal{M}(\mu_1,\dots,\mu_n):= \bigg\{\Q\in \Pi(\mu_1,\cdots,\mu_n)~\bigg|~&\int_{\R^n} H(x_1,\dots,x_i)(x_{i+1}-x_i) \D \Q(x_1,\cdots,x_n)=0 \\
&\text{ for all continuous bounded } H: \R^i \to\R, \text{ for }i = 1,\dots, n-1 \bigg\}
\end{align*} }
denotes the set of martingale measures on $\R^n$ with fixed marginals $\mu_1, \dots, \mu_n${, and where in the above definition $\Pi(\mu_1,\cdots,\mu_n) \subset \mathcal{M}_1(\R^n)$ is defined as the set of all $n$-dimensional Borel probability measures with fixed marginals $\mu_1,\dots,\mu_n$, see for more details, e.g. \cite{villani2009optimal}.}

The value of the MOT problem can either be determined by using the primal formulation outlined in  \eqref{eq_min_MOT} or by maximizing a dual expression of the form
\begin{equation}\label{eq_dual}
\sup_{u_1,\dots,u_n} \sum_{i=1}^n \int _{\R} u_i(x_i) \D \mu_i(x_i)
\end{equation}
 among a certain class of functions $u_i$, compare also \cite{aksamit2019robust, beiglbock2019dual,beiglbock2017complete, cheridito2021martingale, gozlan2017kantorovich,zaev2015monge}. 

The seminal paper \cite{beiglbock2013model} outlines two dual approaches to compute the value of the MOT problem. While the duality provided in \cite[Theorem 1.1]{beiglbock2013model} is formulated and proven in the multi-period case, the convex duality from \cite[Proposition 4.4]{beiglbock2013model} is formulated in the single-period setting, i.e., in the case $n =2$.

In this paper, we reconsider the above mentioned convex duality result for the MOT problem and generalize it to a multi-period formulation.

\subsection{The one-period convex duality result}

{
In the classical optimal transport problem with two marginals (without martingale constraint), the dual problem consists in maximizing \eqref{eq_dual} among functions $u_1,u_2$ fulfilling the inequality constraint $u_1(x_1)+u_2(x_2)\leq c(x_1,x_2)$ for all $(x_1,x_2) \in \R^2$. This turns out to be equivalent to maximizing over pairs of functions $(v_1,v_2)$ where $v_1$ is defined as the \emph{c-transform} or \emph{c-conjugate} of $v_2$, i.e., $v_1(x_1):=\inf_{x_2 \in \R} \{c(x_1,x_2)-v_2(x_2)\}$, and where it is not required to explicilty respect inequality constraints, see also \cite[Chapter 5]{villani2009optimal} for more details.
}

{ It was then shown in \cite[Proposition 4.4]{beiglbock2013model} that for the two marginal MOT problem a similar modified dual formulation can be established which} relies on the notion of { convex \emph{biconjugates}} of real valued functions\footnote{{ Note that Definition~\ref{def_convex_envelope} differs slightly from the definition used in \cite{beiglbock2013model}, where instead of the convex biconjugate the \emph{convex envelope} is considered which is defined identically with the only difference not required to be lower semi-continuous. We additionally require lower semi-continuity to obtain a definition coinciding with the twice applied convex conjugate of a function  as implied by the notation $f^{**}$ which is also used in \cite{beiglbock2013model}. The results of both Proposition~\ref{prop_dual} and \cite[Proposition 4.4]{beiglbock2013model} are valid for either of the definitions, though. }}.

\begin{defn}[{Convex Biconjugate}]\label{def_convex_envelope}
Let $f:\R \rightarrow (-\infty,\infty]$ be some arbitrary function, then we define the  { convex biconjugate }$f^{**}:\R \rightarrow (-\infty,\infty]$ as the largest {lower semi-continuous\footnote{{ A function $f:\R \rightarrow \R$ is called lower semi-continuous  if $\liminf\limits_{x\to x_0} \geq f(x_0)$ for all $x_0 \in \R$. }}  and }convex function smaller or equal to $f$, i.e., it holds for $x \in \R$ that
\[
f^{**}(x):=\sup\{h(x)~|~h:\R \rightarrow \R \text{ {lower semi-continuous and}  convex}, ~h(y) \leq f(y) \text{ for all } y \in \R\}.
\]
\end{defn}
To adapt the notation from \cite{beiglbock2013model},  we define $\R^n \ni (x_1,\dots,x_n) \mapsto S_i(x_1,\dots,x_n):= x_i$ as the canonical process on $\R^n$ and consider first a one-step stochastic process $(S_1,S_2)$ with $S_1 \sim \mu_1,~S_2 \sim\mu_2$. For convenience of the reader we provide below the two-marginal version of the convex duality established in \cite{beiglbock2013model}.

\begin{prop}[\cite{beiglbock2013model}, Proposition 4.4]\label{dual_conv}
Let $c:\R^2 \rightarrow (-\infty,\infty]$ be lower semi-continuous and $c(x,y) \geq -K(1+|x|+|y|) $ for all $x,y \in \R$ and some $K \in \R$. Additionally assume there exists some probability measure $\Q' \in \mathcal{M}(\mu_1,\mu_2)$ such that $\E_\Q'[c(S_1,S_2)] < \infty$. Then, we have
\begin{align*}
&\inf_{\Q \in \mathcal{M}(\mu_1,\mu_2)} \mathbb{E}_\Q[c(S_1,S_2)]=\sup_{u \in L^1(\mu_2) }  \E_{\mu_1}\left[(c(S_1,\cdotb )-u(\cdotb ))^ {**}(S_1)\right]+\E_{\mu_2}[u(S_2)].
\end{align*}
\end{prop}

\section{Main result}
With our main result, formulated in Proposition~\ref{prop_dual}, we extend Proposition~\ref{dual_conv} to the multi-period case with $n>2$, $n \in \N$ marginals. To this end, we consider the $n$-step martingale $S_1,\dots,S_n$ with fixed marginal distributions $S_i \sim \mu_i$ for all $i=1,\dots,n$.

\begin{prop} \label{prop_dual}

Let $c:\R^n \rightarrow (-\infty,\infty]$ be lower semi-continuous and $c(x_1,\dots,x_n) \geq -K(1+|x_1|+\dots+|x_n|) $ for all $x_1,\dots,x_n \in \R$ and some $K \in \R$. Additionally assume there exists some  $\Q' \in \mathcal{M}(\mu_1,\dots,\mu_n)$ such that $\E_{\Q'}[c(S_1,\dots,S_n)] < \infty$. 
We set
\begin{align*}
 c_n:~L^1(\mu_{2})
\times \dots \times L^{1}(\mu_n) \times \R^n &\rightarrow \R \\
 (u_{2},\dots,u_n,x_1,\dots,x_n) &\mapsto c(x_1,\dots,x_n)-\sum_{i=2}^n u_{i}(x_i),
\end{align*}
and 
define inductively for $i = n-1,\dots,1$
\begin{align*}
 c_i:~L^1(\mu_{2})
\times \dots \times L^{1}(\mu_n) \times \R^i &\rightarrow \R \\
 (u_{2},\dots,u_n,x_1,\dots,x_i) &\mapsto \bigg(c_{i+1}(u_{2},\dots,u_n,x_1,\dots,x_i,\cdotb)\bigg)^{**}(x_i).
\end{align*}
Then, we have

\begin{align*}
\inf_{\Q \in \mathcal{M}(\mu_1,\dots,\mu_n)}\mathbb{E}_\Q[c(S_1,S_2,\dots,S_n)] &= \sup_{u_i \in L^1(\mu_i), \atop i = 2,\dots,n}  \E_{\mu_1}\left[c_1(u_2,\dots,u_n,S_1)\right]+\sum_{i=2}^n\E_{\mu_i}[u_i(S_i)].
\end{align*}
\end{prop}
The proof of Proposition~\ref{prop_dual} is provided in Section~\ref{sec_proof}.

\begin{rem}~
\begin{itemize}
\item[(a)] { If $c:\R^n \rightarrow (-\infty,\infty]$ is upper semi-continuous and  fulfills the linear growth condition $c(x_1,\dots,x_n) \leq K(1+|x_1|+\dots+|x_n|)$, then, by applying Proposition~\ref{prop_dual} to $-c$, one obtains a dual expression for the upper bound $\sup_{\Q \in \mathcal{M}(\mu_1,\dots,\mu_n)}\mathbb{E}_\Q[c(S_1,S_2,\dots,S_n)]$ in dependence of the concave biconjugate.\footnote{{ The concave biconjugate denotes the smallest upper semi-continuous concave function larger than $f$.}}}
%
\item[(b)]  {Note that, alternatively we can set \footnote{This means in particular that $c_n$ does not depend on any function $u_i$ for $i=2,\dots,n$.}
$
\widetilde{c}_n \equiv c
$, and 
define recursively for $i = n-1,\dots,1$
\begin{align*}
 \widetilde{c}_i:~L^1(\mu_{i+1})
\times \dots \times L^{1}(\mu_n) \times \R^i &\rightarrow \R \\
 (u_{i+1},\dots,u_n,x_1,\dots,x_i) &\mapsto \bigg(\widetilde{c}_{i+1}(u_{i+2},\dots,u_n,x_1,\dots,x_i,\cdotb)-u_{i+1}(\cdotb)\bigg)^{**}(x_i).
\end{align*}
Then, it follows by induction for $i = n,\dots,2$ that $c_i= \widetilde{c}_i -\sum_{j=2}^i u_j$ and hence $c_1 = \widetilde{c_1}$. Thus, using $(\widetilde{c}_i)_{i=1,\dots,n}$ provides an alternative recursive scheme to construct a dual maximizer.}
{
\item[(c)]
In contrast to the existence of a primal minimizer $\Q^*\in \mathcal{M}(\mu_1,\dots,\mu_n)$, the existence of dual maximizers of the MOT-problem  is not guaranteed without requiring additional assumptions on $c$ (see \cite[Section 4.3]{beiglbock2013model} and \cite{beiglbock2019dual}) or considering a quasi-sure formulation as in \cite{beiglbock2017complete}. If however the dual value from  \cite[Theorem 1.1]{beiglbock2013model} is attained, i.e., if there exist integrable functions $u_i \in L^1(\mu_i)$ for $i=1,\dots,n$ and continuous bounded functions $\Delta_1, \dots, \Delta_{n-1} \in \mathcal{C}_b(\R)$ fulfilling both
\eqref{eq_duality} and $\mathbb{E}_{\Q^*}[c(S_1,\dots,S_n)] =  \sum_{i=1}^n\E_{\mu_i}[u_i(S_i)]$, then we have 
$\sum_{i=1}^n u_i(S_i) + \sum_{j=1}^{n-1} \Delta_j(S_1,\dots,S_j)(S_{j+1}-S_j)= c(S_1,\dots,S_n)$ $\Q^*$-a.s. and hence in the proof of Proposition~\ref{prop_dual}, the inequalities \eqref{eq_duality_1} and \eqref{eq_duality_2} are equalities in the $\Q^*$-a.s sense, implying that $u_1(S_1)=c_1(u_2,\dots,u_n,S_1)$ $\Q^*$-a.s, and therefore that the dual value from Proposition \ref{prop_dual} is attained as well. }
\item[(d)]
Note that for any $\Q \in \mathcal{M}(\mu_1,\dots,\mu_n)$ and for any $u_i \in L^1(\mu_i)$, $i=2,\dots,n$ we have
\begin{equation}\label{eq_cond_ineq}
\E_{\Q}\left[c_1(u_2,\dots,u_n,S_1)+\sum_{i=2}^n u_i(S_i) \middle|~S_1 \right] \leq \E_{\Q}\left[c(S_1,\dots,S_n) \middle|~S_1 \right],
\end{equation}
i.e, the strategy computed in Proposition~\ref{prop_dual} is, conditional on $S_1$, a $\Q$-a.s.\,sub-hedging strategy of $c$. Indeed, by using the martingale property and Jensen's inequality we have
\begin{align*}
&{\E_\Q\bigg[}c_1(u_2,\dots,u_n,S_1)+\sum_{i=2}^n u_i(S_i){~\bigg|~S_1\bigg]}\\
 &= {\E_\Q\bigg[}\left(c_2(u_2,\dots,u_n,S_1, \cdotb)\right)^{**}(\E_{\Q}[S_2|S_1])+\sum_{i=2}^n u_i(S_i){~\bigg|~S_1\bigg]}\\
&\leq \E_\Q\left[ \left(c_2(u_2,\dots,u_n,S_1, \cdotb)\right)^{**}(S_2)+\sum_{i=2}^n u_i(S_i) ~\middle|~S_1\right]
 \\
&\leq \E_\Q\left[ c_2(u_2,\dots,u_n,S_1,S_2)+\sum_{i=2}^n u_i(S_i) ~\middle|~S_1\right]\\
&\leq \E_\Q\left[  \E_\Q\left[c_3(u_2,\dots,u_n,S_1,S_2,S_3)+\sum_{i=2}^n u_i(S_i) ~\bigg|~S_1,S_2\right]~\bigg|~S_1\right] \\
&=\E_\Q\left[  c_3(u_2,\dots,u_n,S_1,S_2,S_3)+\sum_{i=2}^n u_i(S_i) ~\bigg|~S_1\right]\\
\leq \cdots &\leq\E_\Q\left[  c(S_1,\dots,S_n)-\sum_{i=2}^n u_i(S_i) +\sum_{i=2}^n u_i(S_i) ~\bigg|~S_1\right]= \E_\Q\left[  c(S_1,\dots,S_n)~\middle|~S_1\right].
\end{align*}
\end{itemize}
\end{rem}
\subsection{Discussion of Proposition~\ref{prop_dual} and relation to numerics}

As discussed in \cite[Section 4.4]{beiglbock2013model}, the importance and relevance of the c-convex duality theorem stems from the fact that the result allows to construct dual strategies without the need to take inequality constraints into account as it is the case for the sub-hedging duality provided in
\cite[Theorem 1.1]{beiglbock2013model}.
Indeed, the inequality constraints associated to semi-static sub-hedging strategies impose an obstacle when solving the dual formulation of the MOT problem numerically. Therefore, to compute sub-hedging strategies, one usually either considers discretized formulations (\cite{guo2019computational,henry2013automated}) and solves the problem by means of linear programming, or one penalizes any violation of the inequality constraints when parameterizing the functions from the dual expression by appropriate function classes, such as neural networks (\cite{eckstein2021robust, eckstein2021computation}).

%

The simplified construction of hedging strategies based on Proposition~\ref{prop_dual} might open the door to enrich the literature on the numerics of MOT problems by constructing new numerical algorithms.

To showcase one possible way of using the convex duality from Proposition~\ref{prop_dual} numerically, we provide under \href{https://github.com/juliansester/C-Convex}{https://github.com/juliansester/C-Convex} an implementation of an example\footnote{A similar example can be found in \cite[Section 6.4]{alfonsi2019sampling} and \cite[Example 4.3]{sester2020robust}.} from \cite[Section 5.2]{eckstein2021martingale}, where we consider three log-normally marginal distributions. 
The implementation uses Proposition~\ref{prop_dual}, and parameterizes the optimized functions $u_i \in L^1(\mu_i)$, $i=2,3$ by deep neural networks. The objective is then optimized by using the \emph{Adam}-optimizer (\cite{kingma2014adam}) in the \emph{Tensorflow} framework (\cite{abadi2016tensorflow}). To compute the {convex biconjugate} of functions $f$ we make use of the representation of the {convex biconjugate} in terms of the {convex conjugate\footnote{The convex conjugate of a function $f:\R\rightarrow \R$ is defined by $f^*(x)=\sup_{y \in \R} \left\{ y \cdot x -f(y)\right\}$.} $f^*$ applied twice}, given by
$$
\R \ni x \mapsto f^{**}(x)  {=(f^*)^*(x)}= \sup_{m \in \R} \left\{ m \cdot x-\sup_{y \in \R} \left\{ y \cdot m -f(y)\right\} \right\}.
$$
In each iteration we then optimize $\frac{1}{N} \left\{ \sum_{i=1}^N c_1(u_2,u_3,S_1^{(i)}) +\sum_{i=1}^Nu_2(S_2^{(i)}) +\sum_{i=1}^N u_3(S_3^{(i)}) \right\} $ w.r.t.\, the parameters of the deep neural networks representing $u_2,u_3$, where $S_i^{j} \sim \mu_i$ for $i=1,2,3$, $j =1,\dots,N$, and for $N\in \N$ being the batch size. 

{ Another promising approach might be to solve partial differential equations (PDEs) to compute the convex envelope/biconjugate appearing on the dual side of Proposition~\ref{prop_dual}. To this end, one might follow the route of \cite[Corollary 2.2]{henry2017model} which adapts in the two marginal case the optimal transport result from \cite[Proposition 5.48]{villani2003topics} to obtain a link between the solution to the MOT problem and the Hamilton--Jacobi--Bellman equation (which is a non-linear PDE) using the fact that the convex envelope appears naturally as a viscosity solution to certain stochastic control problems which then also satisfies the Hamilton--Jacobi--Bellman PDE, see \cite{pham2005some} and \cite{soner2002stochastic}.
We leave a possible adaption of \cite[Corollary 2.2]{henry2017model} to the $n$-marginal case for future research.}
\section{Proof of Proposition~\ref{prop_dual}} \label{sec_proof}
By using the definition of the { convex biconjugate} we obtain
\begin{equation}\label{eq_smaller_than_1}
\begin{aligned}
&\inf_{\Q \in \mathcal{M}(\mu_1,\dots,\mu_n)}\mathbb{E}_\Q[c(S_1,\dots,S_n)] \\
&=  \inf_{\Q \in \mathcal{M}(\mu_1,\dots,\mu_n)} \sup_{u_i \in L^1(\mu_i), \atop i = 2,\dots,n}  \mathbb{E}_\Q \left[c(S_1,\dots,S_{n-1},S_n)-\sum_{i=2}^nu_i(S_i)\right] +\sum_{i=2}^n\E_{\mu_i}[u_i(S_i)]  \\
&\geq  \inf_{\Q \in \mathcal{M}(\mu_1,\dots,\mu_n)} \sup_{u_i \in L^1(\mu_i), \atop i = 2,\dots,n} \mathbb{E}_\Q\left[\bigg(c(S_1,\dots,S_{n-1},\cdotb)-\sum_{i=2}^{n-1}u_i(S_i)-u_n(\cdotb)\bigg)^{**}(S_n)\right]  +\sum_{i=2}^n\E_{\mu_i}[u_i(S_i)]
\end{aligned}
\end{equation}
By \eqref{eq_smaller_than_1}, by the tower property of the conditional expectation and by Jensen's inequality we have
\begin{equation}\label{eq_smaller_than_2}
\begin{aligned}
&\inf_{\Q \in \mathcal{M}(\mu_1,\dots,\mu_n)}\mathbb{E}_\Q[c(S_1,\dots,S_n)]\\  &\geq  \inf_{\Q \in \mathcal{M}(\mu_1,\dots,\mu_n)} \sup_{u_i \in L^1(\mu_i), \atop i = 2,\dots,n}  \mathbb{E}_\Q\left[\mathbb{E}_\Q\left[\bigg(c(S_1,\dots,S_{n-1},\cdotb)-\sum_{i=2}^{n-1}u_i(S_i)-u_n(\cdotb)\bigg)^{**}(S_n) ~\middle|~S_{1},\dots,S_{n-1}\right]\right] \\
&\hspace{5cm}+\sum_{i=2}^n\E_{\mu_i}[u_i(S_i)] \\
&\geq   \inf_{\Q \in \mathcal{M}(\mu_1,\dots,\mu_n)}\sup_{u_i \in L^1(\mu_i), \atop i = 2,\dots,n} \mathbb{E}_\Q\left[ \bigg(c(S_1,\dots,S_{n-1},\cdotb)-\sum_{i=2}^{n-1}u_i(S_i)-u_n(\cdotb)\bigg)^{**}(\E_{\Q} \left[S_n~\middle|~S_{1},\dots,S_{n-1}\right])\right]   \\
&\hspace{5cm}+\sum_{i=2}^n\E_{\mu_i}[u_i(S_i)].
\end{aligned}
\end{equation}
Due to \eqref{eq_smaller_than_2} and the martingale property of  $\Q \in \mathcal{M}(\mu_1,\dots,\mu_n)$ we obtain
\begin{equation*}
\begin{aligned}
&\inf_{\Q \in \mathcal{M}(\mu_1,\dots,\mu_n)}\mathbb{E}_\Q[c(S_1,\dots,S_n)] \\
&\geq    \inf_{\Q \in \mathcal{M}(\mu_1,\dots,\mu_n)} \sup_{u_i \in L^1(\mu_i), \atop i = 2,\dots,n}\mathbb{E}_\Q\left[ \bigg(c(S_1,\dots,S_{n-1},\cdotb)-\sum_{i=2}^{n-1}u_i(S_i)-u_n(\cdotb)\bigg)^{**}(S_{n-1})\right]  +\sum_{i=2}^n\E_{\mu_i}[u_i(S_i)]  \\
&=   \inf_{\Q \in \mathcal{M}(\mu_1,\dots,\mu_n)} \sup_{u_i \in L^1(\mu_i), \atop i = 2,\dots,n}  \mathbb{E}_\Q\left[ c_{n-1}(u_2,\dots,u_n,S_1,\dots,S_{n-1})\right]   +\sum_{i=2}^n\E_{\mu_i}[u_i(S_i)]  .
\end{aligned}
\end{equation*}
By proceeding with the same arguments we have 
\begin{align*}
&\inf_{\Q \in \mathcal{M}(\mu_1,\dots,\mu_n)}\mathbb{E}_\Q[c(S_1,\dots,S_n)] \\
&\geq   \inf_{\Q \in \mathcal{M}(\mu_1,\dots,\mu_n)} \sup_{u_i \in L^1(\mu_i), \atop i = 2,\dots,n}  \mathbb{E}_\Q\left[ \left(c_{n-1}(u_2,\dots,u_n,S_1,\dots,S_{n-2},\cdotb)\right)^{**}(S_{n-1})\right]   +\sum_{i=2}^n\E_{\mu_i}[u_i(S_i)] \\
&\geq   \inf_{\Q \in \mathcal{M}(\mu_1,\dots,\mu_n)} \sup_{u_i \in L^1(\mu_i), \atop i = 2,\dots,n}  \mathbb{E}_\Q\left[ \left(c_{n-1}(u_2,\dots,u_n,S_1,\dots,S_{n-2},\cdotb)\right)^{**}(S_{n-2})\right]   +\sum_{i=2}^n\E_{\mu_i}[u_i(S_i)] \\
&=   \inf_{\Q \in \mathcal{M}(\mu_1,\dots,\mu_n)} \sup_{u_i \in L^1(\mu_i), \atop i = 2,\dots,n}  \mathbb{E}_\Q\left[ c_{n-2}(u_2,\dots,u_n,S_1,\dots,S_{n-2})\right]   +\sum_{i=2}^n\E_{\mu_i}[u_i(S_i)],
\end{align*}
which implies inductively that 
\begin{equation}\label{eq_new_eq_revision1}
\begin{aligned}
&\inf_{\Q \in \mathcal{M}(\mu_1,\dots,\mu_n)}\mathbb{E}_\Q[c(S_1,\dots,S_n)] \geq \sup_{u_i \in L^1(\mu_i), \atop i = 2,\dots,n}  \E_{\mu_1}\left[c_1(u_2,\dots,u_n,S_1)\right]+\sum_{i=2}^n\E_{\mu_i}[u_i(S_i)].
 \end{aligned}
\end{equation}
Next, note that if there exist integrable functions $u_i \in L^1(\mu_i)$ for $i=1,\dots,n$ and continuous bounded functions $\Delta_1, \dots, \Delta_{n-1} \in \mathcal{C}_b(\R)$ such that 
\begin{equation}\label{eq_duality}
\begin{aligned}
\Psi_{(u_i),(\Delta_j)}(x_1,\dots,x_n):&=\sum_{i=1}^n u_i(x_i) + \sum_{j=1}^{n-1} \Delta_j(x_1,\dots,x_{j})(x_{j+1}-x_j)\leq c(x_1,\dots,x_{n-1},x_n) 
\end{aligned}
\end{equation}
for all $ (x_1,\dots,x_n) \in \R^n$,
then we see that
\begin{equation}\label{eq_duality_1}
\begin{aligned}
&u_1(x_1)+ \sum_{j=1}^{n-1} \Delta_j(x_1,\dots,x_{j})(x_{j+1}-x_j) \\
=&\left(u_1(x_1)+ \sum_{j=1}^{n-2} \Delta_j(x_1,\dots,x_{j})(x_{j+1}-x_j) + \Delta_{n-1}(x_1,\dots,x_{n-1})(\cdotb-x_j) \right)^{**}(x_n) \\
\leq &\left( c(x_1,\dots,x_{n-1},\cdotb)-\sum_{i{=2}}^{n-1} u_i(x_i) -u_n(\cdotb) \right)^{**}(x_n)    
\end{aligned}
\end{equation}
 for all $(x_1,\dots,x_n) \in \R^n$, since the left-hand side of \eqref{eq_duality_1} is linear (and hence convex) in $x_n$. By setting $x_n := x_{n-1}$, we obtain  from \eqref{eq_duality_1} that
 \begin{equation}\label{eq_duality_2}
u_1(x_1) + \sum_{j=1}^{n-2} \Delta_j(x_1,\dots,x_{j})(x_{j+1}-x_j) \leq  c_{n-1}(u_2,\dots,u_n,x_1,\dots,x_{n-1})  
\end{equation}
for all $(x_1,\dots,x_{n-1}) \in \R^{n-1}$. Repeating inductively the argument leading to \eqref{eq_duality_2} shows that \eqref{eq_duality} implies 
 \begin{equation}\label{eq_duality_3}
 u_1(x_1) \leq c_1(u_2,\dots,u_n,x_1) \text{ for all } x_1 \in \R.
 \end{equation}
 Hence, we have {due to \eqref{eq_new_eq_revision1}} and \eqref{eq_duality_3}  that\footnote{The notation $\Psi_{(u_i),(\Delta_j)} \leq c$ means $\Psi_{(u_i),(\Delta_j)}(x) \leq c(x)$ for all $x \in \R^n$. }
 \begin{equation}\label{eq_smaller_than_3}
\begin{aligned}
&\inf_{\Q \in \mathcal{M}(\mu_1,\dots,\mu_n)}\mathbb{E}_\Q[c(S_1,\dots,S_n)] \\
 &\geq \sup_{u_i \in L^1(\mu_i), \atop i = 2,\dots,n}  \E_{\mu_1}\left[c_1(u_2,\dots,u_n,S_1)\right]+\sum_{i=2}^n\E_{\mu_i}[u_i(S_i)]\\
  &\geq \sup_{~{u_i \in L^1(\mu_i), i = 1,\dots,n ,\atop
  \Delta_j \in \mathcal{C}_b(\R), j =1,\dots,n-1 }: \Psi_{(u_i),(\Delta_j)} \leq c}   ~\E_{\mu_1}\left[c_1(u_2,\dots,u_n,S_1)\right]+\sum_{i=2}^n\E_{\mu_i}[u_i(S_i)] \\
&\geq \sup_{~{u_i \in L^1(\mu_i), i = 1,\dots,n ,\atop
  \Delta_j \in \mathcal{C}_b(\R), j =1,\dots,n-1 }: \Psi_{(u_i),(\Delta_j)} \leq c}   ~\sum_{i=1}^n\E_{\mu_i}[u_i(S_i)] = \inf_{\Q \in \mathcal{M}(\mu_1,\dots,\mu_n)}\mathbb{E}_\Q[c(S_1,\dots,S_n)],
 \end{aligned}
\end{equation}
where in the last step we used the duality result from \cite[Theorem 1.1]{beiglbock2013model}. The applicability of \cite[Theorem 1.1]{beiglbock2013model} is ensured by the linear growth condition and the lower semi-continuity imposed on $c$. Hence, all inequalities in \eqref{eq_smaller_than_3} are indeed equalities and the assertion follows.
\bibliography{literature}

\bibliographystyle{plain}
\end{document}